\documentclass{article}
\usepackage[latin1]{inputenc}
\usepackage[T1]{fontenc}
\usepackage{amssymb}
\usepackage{showlabels}
\usepackage{times}

\usepackage{color}
\definecolor{DarkOlive}{rgb}{0.1047,0.2412,0.0064}
\definecolor{FireBrick}{rgb}{0.5812,0.0074,0.0083}
\definecolor{RoyalBlue}{rgb}{0.0236,0.0894,0.6179}
\definecolor{RoyalGreen}{rgb}{0.0236,0.6179,0.0894}
\definecolor{RoyalRed}{rgb}{0.6179,0.0236,0.0894}
\definecolor{LightBlue}{rgb}{0.8544,0.9511,1.0000}
\definecolor{Black}{rgb}{0.0,0.0,0.0}
\definecolor{MidnightBlue}{rgb}{0.0035,0.0020,0.1363}
\definecolor{FireBrick3}{rgb}{0.5812,0.0074,0.0083}
\definecolor{FireBrick4}{rgb}{0.2156,0.0023,0.0035}
\definecolor{Blue2}{rgb}{0.0000,0.0000,0.8644}
\definecolor{Navy}{rgb}{0.0000,0.0000,0.1927}
\definecolor{MediumBlue}{rgb}{0.0000,0.0000,0.6179}
\usepackage[
        a4paper=true,bookmarks=false,pdftitle={Generic Computations in
        Finite Groups of Lie Type},
        colorlinks=true,backref=false,breaklinks=true,linkcolor=MediumBlue,
        citecolor=FireBrick3,filecolor=RoyalRed,
        urlcolor=Blue2,pagecolor=MediumBlue]{hyperref}

\newtheorem{Thm}{Theorem}

\newtheorem{Pro}[Thm]{Proposition}
\newtheorem{La}[Thm]{Lemma}
\newtheorem{Rem}[Thm]{Remark}

\newenvironment{Prf}{\noindent\textbf{Proof.}}{\hfill $\Box$ \medskip}

\newcommand{\Z}{\mathbb{Z}}

\begin{document}
\title{Elements with only negative cycles in Weyl groups of type $B$ and $D$}
\author{Frank Lübeck}
\maketitle
\begin{abstract}
We  give  a  precise  formula  and a  simple  asymptotic  function  for  the
proportion of elements with only negative  (or only positive) cycles in Weyl
groups of type $B,C$ and $D$.
\end{abstract}

For $n \in  \Z_{\geq 1}$ let $W(B_n)$  and $W(D_n)$ be the  Weyl groups with
$n$ Coxeter  generators of  type $B$  and $D$,  respectively. $W(B_n)$  is a
wreath product of a cyclic group of  order $2$ with a symmetric group on $n$
points, and $W(D_n)$ is a normal subgroup of index $2$ in $W(B_n)$.

We      recall      some      well       known      facts,      see      for
example~\cite[7.]{CarterConjugacyWeyl}. The  group  $W(B_n)$ has  a  natural
representation   as  permutation   group  on   $2n$  points   $\{1,2,\ldots,
n,1',\ldots,n'\}$; the  group is  generated by the  transpositions $(i,i')$,
$1\leq  i\leq n$,  and permutations  of form  $\pi\pi'$, where  $\pi$ is  an
element  of the  symmetric  group on  $\{1,\ldots,n\}$  and $\pi'$  permutes
$\{1',\ldots,n'\}$ with $\pi'(i')  = \pi(i)'$. We consider the  cycles of an
element $w \in W(B_n)$.  If such a cycle contains a point  $i$ and $i'$ then
the cycle has even length $2k$ and we call this cycle a \emph{negative cycle
of length $k$}. If  a cycle of length $k$ contains a point  $i$ and not $i'$
then there is  another cycle of length $k$ containing  $i'$. Such cycles are
called \emph{positive cycles of length $k$}.

The normal subgroup $W(D_n)$ consists of the elements with an even number of
negative cycles.

There is  a natural  surjection $\sigma:  W(B_n) \to  S_n$ to  the symmetric
group on the $n$ pairs $\{i,i'\}$ (each $w \in W(B_n)$ induces a permutation
on these pairs).

\begin{La}
For $x \in S_n$ we consider  the $2^n$ preimages $M_x = \sigma^{-1}(x)$. Let
$x=x_1\cdots x_k$ be  the decomposition of $x$ into disjoint  cycles. For $w
\in M_x$ let  $s(w) \in \{1,-1\}^k$ be  the sequence of $k$  signs such that
the $j$-th entry is $1$ if two positive  cycles of $w$ map to $x_j$ and $-1$
if a negative cycle of $w$ maps to $x_j$. For $s \in \{1,-1\}^k$ let $M_x(s)
= \{w\in  M_x\mid\; s(w) =  s\}$. Then all the  sets $M_x(s)$ have  the same
cardinality $2^{n-k}$.
\end{La}
\begin{Prf}
Let $s_0  = (1,1,\ldots,1)$ and  $s \in \{1,-1\}^k$  such that it  has entry
$-1$  in positions  $j_1,\ldots,j_l$. For  each  $j_m$, $1\leq  m\leq l$  we
choose a pair  $\{i_m,i_m'\}$ in the cycle $x_{j_m}$. Then  the map $M_x \to
M_x$, $w \mapsto (i_1,i_1')\cdots(i_l,i_l')w$, is a bijection of $M_x$ which
maps $M_x(s_0)$ to $M_x(s)$.
\end{Prf}

\begin{Pro}
(a) The proportion $p(n)$ of elements in $W(B_n)$ with only negative cycles is 
\[ p(n) = \frac{1\cdot 3\cdot 5 \cdots (2n-1)}{2\cdot 4\cdot 6 \cdots
(2n)} = \frac{(2n-1)!!}{|W(B_n)|}.\]

(b) The proportion $p^+(n)$ of elements of $W(D_n)$ with only negative cycles is
\[ p^+(n) = p(n) \cdot \frac{2n-2}{2n-1}. \]

(c) The proportion $p^-(n)$ of elements in the coset 
$W(B_n)\setminus W(D_n)$ with  only negative cycles is 
\[ p^-(n) = p(n) \cdot \frac{2n}{2n-1}. \]
\end{Pro}
\begin{Prf}
(a)
The number $S_1(n,k)$  of permutations in $S_n$ with $k$  cycles is given by
the unsigned  Stirling numbers of the  first kind. These can  be obtained as
coefficients of a polynomial as follows 
(see, e.g.,~\cite[III.2.4]{AignerKombinatorik})
\[ \chi_n(X) := \sum_{k=0}^n S_1(n,k) X^k = X (X+1)\cdots (X+n-1). \]

From the lemma  we know that for  each element of $S_n$ with  $k$ cycles its
preimage in $W(B_n)$ contains $2^{n-k}$  elements with only negative cycles.
So, all together we have
\[ \sum_{k=0}^n S_1(n,k) \cdot 2^{n-k} = 2^n \cdot \chi_n(1/2) = (2n-1)!! \]
elements with only negative cycles in  $W(B_n)$. Since $|W(B_n)| = 2^n \cdot
n!$ we get~(a).

(b) 
To find the  elements with only negative cycles in  the subgroup $W(D_n)$ we
can proceed as for part~(a) but only sum over the terms with even $k$. To do
so, we use
\[ \sum_{k=0,\, k \textrm{\scriptsize\  even}}^n S_1(n,k) \cdot 2^{n-k}
= 2^{n-1} (\chi_n(1/2) + \chi_n(-1/2)). \]
Using again the polynomial for $\chi_n(X)$ we get the formula for $p^+(n)$.

(c) This is analogous to~(b), now using 
\[ \sum_{k=0,\, k \textrm{\scriptsize\  odd}}^n S_1(n,k) \cdot 2^{n-k}
= 2^{n-1} (\chi_n(1/2) - \chi_n(-1/2)). \]
\end{Prf}

\begin{Rem}
Let 
\[ h(n) = \frac{\left(1+\frac{1}{22 n}\right)}{\sqrt{\pi n}}.\]
Then we have for  all $n \in \Z_{\geq 1}$ that $p(n)  < h(n)$ and $p(n)/h(n)
\to 1$ for $n \to \infty$.
\end{Rem}
\begin{Prf}
This follows from writing 
\[ p(n) = \frac{(2n-1)!!}{2^n n!} = \frac{(2n)!}{(2^n n!)^2} \]
and the well known Stirling approximation of factorials:
\[ \sqrt{2\pi n} \left(\frac{n}{e}\right)^n < n! < \left(1+\frac{1}{11\cdot
n}\right) \sqrt{2\pi n} \left(\frac{n}{e}\right)^n.\]
\end{Prf}

\emph{Remark.} This short note answers a question by Gerhard Hiss. The above
results will be interpreted as proportions of certain elements in finite 
groups of Lie type in~\cite[7.6]{HissHusenMagaard}.


\begin{thebibliography}{Car72}

\bibitem[Aig75]{AignerKombinatorik}
M.~Aigner.
\newblock {\em Kombinatorik I}.
\newblock Springer-Verlag, Berlin, 1975.

\bibitem[Car72]{CarterConjugacyWeyl}
R.W. Carter.
\newblock Conjugacy classes in the {W}eyl group.
\newblock {\em Compositio Math.}, 25:1--59, 1972.

\bibitem[HHM]{HissHusenMagaard}
G.~Hiss, W.~J. Husen, and K.~Magaard.
\newblock Imprimitive irreducible modules for finite quasisimple groups.
\newblock {\em preprint}.

\end{thebibliography}

\end{document}